\documentclass[11pt]{amsart}

\usepackage[T1]{fontenc}
\usepackage{amssymb}
\usepackage[margin=1.2in]{geometry}
\usepackage{microtype}

\newtheorem{theorem}{Theorem}
\newtheorem{proposition}[theorem]{Proposition}
\newtheorem{lemma}[theorem]{Lemma}
\theoremstyle{remark}

\newcommand{\Area}{\operatorname{Area}}
\newcommand{\diam}{\operatorname{diam}}

\begin{document}

\title[An improved isodiametric inequality]{An improved isodiametric inequality for convex surfaces}

\author{Mat\v{e}j Stehl\'{\i}k}
\address{Universit\'e Paris Cit\'e, CNRS, IRIF, F-75013 Paris, France}
\email{matej@irif.fr}


\begin{abstract}
  Alexandrov conjectured that every closed convex surface of area $A$ and intrinsic diameter $D$ satisfies $A/D^2\leq\pi/2$.  Using an averaging argument over the curvature measure, we prove that every closed convex surface satisfies $A/D^2\leq2\pi/3$, improving the best previously known bound.
\end{abstract}

\maketitle

\section{Introduction}

A \emph{closed convex surface} is the boundary of a convex body in $\mathbb R^3$, or a doubly covered compact convex plane domain (a degenerate convex surface), equipped with its intrinsic metric. Alexandrov~\cite[p.~417]{Alexandrov} conjectured that every closed convex surface satisfies
\begin{equation}\label{eq:conj}
  \frac{A}{D^2}\leq\frac{\pi}{2}\approx1.571,
\end{equation}
where $A$ is its area and $D$ its intrinsic diameter. Equality holds for the doubly covered flat disc. The conjecture is often stated, as in~\cite{FreitasKrejcirik,Shioya}, for smooth closed Riemannian two-spheres of nonnegative Gaussian curvature; the doubly covered disc is then a non-smooth limit of such surfaces, and the constant $\pi/2$ is approached but not attained. The conjecture is wide open. It is known for several restricted classes, including surfaces of revolution and ellipsoids; see~\cite{FreitasKrejcirik}.

Bishop's inequality, applied at any point, gives $A/D^2\leq\pi$. Area--diameter estimates for nonnegatively curved two-spheres and projective planes, and partial results on the conjecture, were obtained by Sakai~\cite{Sakai}, Shioya~\cite{Shioya93} and Zalgaller~\cite{Zalgaller}. For general smooth convex surfaces, Calabi and Cao~\cite{CalabiCao} used eigenvalue estimates to prove $A/D^2\leq8/\pi\approx2.546$, and Shioya~\cite{Shioya} proved $A/D^2\leq(\sqrt{21}-3)\pi/2\approx2.486$, which to our knowledge is the best previously known bound. Our main result is the following improvement.

\begin{theorem}\label{thm:main}
  Every closed convex surface with area $A$ and intrinsic diameter $D$ satisfies
  \begin{equation}\label{eq:main}
    \frac{A}{D^2}\leq\frac{2\pi}{3}\approx2.094.
  \end{equation}
\end{theorem}

The heart of the paper is the corresponding result for smooth surfaces, in the setting of~\cite{FreitasKrejcirik,Shioya}.

\begin{theorem}\label{thm:smooth}
  Let $M$ be a smooth closed Riemannian surface homeomorphic to the two-sphere, with nonnegative Gaussian curvature, area $A$ and diameter $D$. Then $A/D^2<2\pi/3$.
\end{theorem}

We prove Theorem~\ref{thm:smooth} in Section~\ref{sec:proof}, and deduce Theorem~\ref{thm:main} from it by approximation at the end of that section. The ingredients are the same as in Shioya's argument: the Fiala--Hartman inequality for geodesic circles, Gauss--Bonnet and Bishop--Gromov comparison. The difference lies in the order in which they are combined; see Section~\ref{sec:shioya}. In Section~\ref{sec:projective} we note that the same argument also improves Shioya's bound for the real projective plane.

\section{The proof}\label{sec:proof}

From now until the end of the proof of Theorem~\ref{thm:smooth}, $M$ is as in Theorem~\ref{thm:smooth}, $d$ is its distance function and $K \geq 0$ its Gaussian curvature. By the Gauss--Bonnet theorem,
\begin{equation}\label{eq:GB}
  \int_MK\,dA=4\pi.
\end{equation}
For $p\in M$ and $r\geq0$ let $B(p,r)=\{q\in M:d(p,q)<r\}$ be the open ball, and let
\[
  R_p=\max_{z\in M}d(p,z)
\]
be the eccentricity of $p$, so that $R_p\leq D$. For $r>0$ the boundary $\partial B(p,r)$ is the distance circle $\{q:d(p,q)=r\}$. Following~\cite[Section~2.1]{BerardCastillon}, let $L_p$ be the function which, outside a closed set of measure zero, equals the length of $\partial B(p,r)$, extended to a right-continuous function on $[0,\infty)$ with $L_p(0)=0$; its jumps are downward. (The setting of~\cite{BerardCastillon} is that of complete open surfaces; the proof of Proposition~\ref{prop:FH} explains how it applies here.) Since $L_p(r)$ agrees with the length of $\partial B(p,r)$ for almost every $r$, and $|\nabla d(p,\cdot)|=1$ almost everywhere, the coarea formula gives $\Area(B(p,t))=\int_0^tL_p(s)\,ds$.

\begin{proposition}[Fiala, Hartman]\label{prop:FH}
  For every $p\in M$ and every $0<t<R_p$,
  \begin{equation}\label{eq:FH}
    L_p(t)\leq2\pi t-\int_0^t\!\int_{B(p,s)}K\,dA\,ds.
  \end{equation}
\end{proposition}

\begin{proof}
  For complete, connected, orientable surfaces, the integrated Fiala--Hartman inequality~\cite{Fiala}, \cite[Proposition~6.1]{Hartman}, \cite[Proposition~3.7]{ShiohamaTanaka}, in the form recorded in~\cite[(2.2), (2.8)]{BerardCastillon}, states that for $0\leq a<b$,
  \[
    L_p(b)-L_p(a)\leq2\pi(b-a)\sup_{a\leq s\leq b}\chi\bigl(B(p,s)\bigr)
      -\int_a^b\!\int_{B(p,s)}K\,dA\,ds,
  \]
  where $\chi$ is the Euler characteristic, with the convention that it equals $-\infty$ for a surface of infinite topology~\cite[Section~2.1]{BerardCastillon}; see also~\cite[Theorem~4.4.1]{SST}.

  We first reduce to this setting. Fix $0<t<R_p$, choose $z\in M$ with $d(p,z)>t$, and choose $t'$ with $t<t'<d(p,z)$ and a neighbourhood $U$ of $z$ whose closure is disjoint from the closed ball $\overline{B}(p,t')$.  The surface $M\setminus\{z\}$ is homeomorphic to the plane. On it, let $g'=\varphi g$, where $\varphi$ is smooth, $\varphi\geq1$, $\varphi=1$ outside $U$, and $\varphi$ grows fast enough near $z$ that $g'$ is complete (for instance $\varphi(x)\geq d(z,x)^{-2}$ in a small normal neighbourhood of $z$). Let $d'$ be the distance of $g'$.  Since $g'\geq g$, we have $d'\geq d$. Conversely, if $d(p,q)<s\leq t'$, a minimizing $g$-geodesic from $p$ to $q$ lies in $B(p,s)$, where $g'=g$, so $d'(p,q)\leq d(p,q)$. Hence the $g'$-balls of radius $s\leq t'$ about $p$ coincide with the balls $B(p,s)$, and the distance circles, their lengths and the curvature on these balls are the same for $g$ and $g'$. Since $t'>t$, this holds on a right neighbourhood of $t$, so the right-continuous representatives of the length functions for $g$ and $g'$ agree on $[0,t]$. The inequality above therefore holds for $M$ and $0\leq a<b\leq t$.

  Take $a=0$ and $b=t$. The ball $B(p,0)$ is empty. For $0<s\leq t$, the ball $B(p,s)$ is connected, since it is a union of minimizing geodesic segments from $p$. It is an open proper subset of the connected surface $M$, since it misses every point at distance $R_p$ from $p$, so it is not closed, hence not compact.  Under the conventions of~\cite[Section~2.1]{BerardCastillon}, a connected noncompact surface has Euler characteristic at most~$1$ (equal to $-\infty$ in the case of infinite topology). Hence the supremum is at most~$1$, and \eqref{eq:FH} follows.
\end{proof}

Since $K\geq0$, the Bishop--Gromov comparison theorem~\cite{Chavel} shows that $r\mapsto\Area(B(q,r))/r^2$ is nonincreasing for every $q\in M$. For $0<r\leq D<r'$ this gives $\Area(B(q,r))/r^2\geq\Area(B(q,r'))/r'^2=A/r'^2$, since $B(q,r')=M$; letting $r'\downarrow D$, we obtain
\begin{equation}\label{eq:BG}
  \frac{\Area(B(q,r))}{A}\geq\frac{r^2}{D^2}, \qquad0\leq r\leq D.
\end{equation}

We shall also use the following layer-cake identity. If $0 \leq a \leq D$, then $(D-a)^2=\int_0^D2(D-r)\,\mathbf 1[a<r]\,dr$. Taking $a=d(p,q)$ and integrating over $p\in M$ gives, by Tonelli's theorem,
\begin{equation}\label{eq:layercake}
  \int_M\bigl(D-d(p,q)\bigr)^2\,dA(p) = \int_0^D2(D-r)\Area(B(q,r))\,dr.
\end{equation}

\begin{lemma}\label{lem:layer}
  For every $p\in M$,
  \begin{equation}\label{eq:layer}
    A\leq\pi R_p^2 - \frac12 \int_MK(q)\bigl(R_p-d(p,q)\bigr)^2\,dA(q).
  \end{equation}
\end{lemma}

\begin{proof}
  The set of points at distance exactly $R_p$ from $p$ has measure zero, so $A$ is the area of $B(p,R_p)$, that is, $A=\int_0^{R_p}L_p(t)\,dt$. Applying \eqref{eq:FH} for $0<t<R_p$,
  \[
   A\leq\int_0^{R_p}\Bigl(2\pi t-\int_0^t\!\int_{B(p,s)}K\,dA\,ds\Bigr)dt.
  \]
  The first term contributes $\pi R_p^2$.  For fixed $q$, the set of $(s,t)$ with $q\in B(p,s)$ and $0\leq s\leq t\leq R_p$ is the triangle $d(p,q)<s\leq t\leq R_p$, of area $\tfrac12(R_p-d(p,q))^2$. Since $K\geq0$, Tonelli's theorem gives \eqref{eq:layer}.
\end{proof}

For $p\in M$, define the \emph{curvature potential} of $p$ as
\[
  P_D(p)=\int_MK(q)\bigl(D-d(p,q)\bigr)^2\,dA(q).
\]

\begin{lemma}\label{lem:average}
  There is a point $p\in M$ such that
  \begin{equation}\label{eq:average}
    P_D(p)\geq\frac{2\pi}{3}D^2.
  \end{equation}
\end{lemma}

\begin{proof}
  Fix $q\in M$.  By \eqref{eq:layercake} and \eqref{eq:BG},
  \begin{equation}\label{eq:avgq}
    \frac1A\int_M\bigl(D-d(p,q)\bigr)^2\,dA(p) = \int_0^D2(D-r)\frac{\Area(B(q,r))}{A}\,dr \geq\int_0^D2(D-r)\frac{r^2}{D^2}\,dr = \frac{D^2}{6}.
  \end{equation}
  Multiply by $K(q) \geq 0$, integrate over $q$, and use Tonelli's theorem and \eqref{eq:GB}:
  \[
    \frac1A\int_MP_D(p)\,dA(p)\geq4\pi\cdot\frac{D^2}{6}=\frac{2\pi}{3}D^2.
  \]
  Some point $p$ therefore satisfies \eqref{eq:average}.
\end{proof}

\begin{proof}[Proof of Theorem~\ref{thm:smooth}]
  Choose $p$ as in Lemma~\ref{lem:average}, and put
  \[
    E=R_p,\qquad\delta=D-E\geq0,\qquad
    U=\int_MK(q)\bigl(E-d(p,q)\bigr)^2\,dA(q).
  \]
  Since $D-d(p,q)=\bigl(E-d(p,q)\bigr)+\delta$, Minkowski's inequality in $L^2(K\,dA)$ and \eqref{eq:GB} give
  \[
    \sqrt{P_D(p)}\leq\sqrt U+2\sqrt\pi\,\delta.
  \]
  Combined with \eqref{eq:average}, this gives
  \begin{equation}\label{eq:U}
    U\geq\bigl(a-b\delta\bigr)_+^2,
    \qquad a=\sqrt{\tfrac{2\pi}{3}}\,D,\quad b=2\sqrt\pi,
  \end{equation}
  where $t_+=\max\{t,0\}$. Lemma~\ref{lem:layer} at $p$ and \eqref{eq:U} imply
  \begin{equation}\label{eq:F}
    A\leq F(\delta):=\pi(D-\delta)^2-\tfrac12\bigl(a-b\delta\bigr)_+^2.
  \end{equation}
  For $0\leq\delta\leq a/b=D/\sqrt6$,
  \[
    F'(\delta) = -2\pi(D-\delta)+b(a-b\delta) = 2\pi D\Bigl(\sqrt{\tfrac23}-1\Bigr)-2\pi\delta<0.
  \]
  For $a/b\leq\delta\leq D$, $F(\delta)=\pi(D-\delta)^2$ is again decreasing. Since $F$ is continuous, it is maximised at $\delta=0$, and
  \[
    A\leq F(0)=\pi D^2-\tfrac12\cdot\tfrac{2\pi}{3}D^2,
    \qquad\text{that is,}\qquad\frac{A}{D^2}\leq\frac{2\pi}{3}.
  \]

  It remains to show that the inequality is strict. Let $q$ be a point with $K(q)>0$. The function $f(r)=\Area(B(q,r))/r^2$ is nonincreasing, and $f(r)=\pi-\tfrac{\pi}{12}K(q)r^2+O(r^3)$ as $r\to0$~\cite{Gray}.  If equality held in \eqref{eq:BG} for almost every $r\in(0,D)$, then $f$ would be constant on $(0,D)$, hence equal to $\pi$, contradicting the expansion. So \eqref{eq:avgq} is strict for every such $q$. Since $\{K>0\}$ has positive measure by \eqref{eq:GB}, integrating as in the proof of Lemma~\ref{lem:average} gives
  \[
    \frac1A\int_MP_D(p)\,dA(p)>\frac{2\pi}{3}D^2,
  \]
  so we may choose $p$ with $P_D(p)>\tfrac{2\pi}{3}D^2$. Then \eqref{eq:U} holds with $a$ replaced by $a'=\sqrt{P_D(p)}>a$, and $A\leq\widetilde F(\delta)$, where $\widetilde F$ is defined as $F$ but with $a'$ in place of $a$; note that $\widetilde F\leq F$. If $\delta=0$, then $A\leq\widetilde F(0)<F(0)$. If $\delta>0$, then $A\leq F(\delta)<F(0)$, since $F$ is strictly decreasing. In both cases $A<F(0)=\tfrac{2\pi}{3}D^2$.
\end{proof}

\begin{proof}[Proof of Theorem~\ref{thm:main}]
  Let $\Sigma$ be a closed convex surface with area $A$ and intrinsic diameter $D$.  Suppose first that $\Sigma$ is degenerate, the double of a compact convex plane domain $\Omega$. Then $A=2\Area(\Omega)$. Intrinsic distances on $\Sigma$ are at least Euclidean distances, so $D\geq\diam(\Omega)$, and the planar isodiametric inequality $\Area(\Omega)\leq\pi\diam(\Omega)^2/4$~\cite{Bieberbach} gives $A/D^2\leq\pi/2<2\pi/3$.

  Suppose now that $\Sigma=\partial C$ for a convex body $C$ with nonempty interior. Let $\varepsilon_n\downarrow0$, let $T_n$ be the smoothing map of~\cite[Theorem~3.4.1]{Schneider} with parameter $\varepsilon_n$, and put $C_n=T_nC+\varepsilon_nB$, where $B$ is the unit ball. After a translation, we may choose $R>0$ with $C\subset RB$; then~\cite[Theorem~3.4.1(c)]{Schneider} gives $\delta(C,C_n)\leq(R+1)\varepsilon_n$, so $C_n\to C$ in the Hausdorff metric $\delta$. The support function of $C_n$ is $h_{T_nC}+\varepsilon_n$, which is smooth on $\mathbb R^3\setminus\{0\}$ with all radii of curvature at least $\varepsilon_n$; hence $\partial C_n$ is a smooth surface of positive Gaussian curvature. Theorem~\ref{thm:smooth} therefore applies to $\partial C_n$. By Alexandrov's convergence theorem~\cite[p.~106]{Alexandrov}, if points $X_n,Y_n\in\partial C_n$ converge to $X,Y\in\Sigma$, then their intrinsic distances on $\partial C_n$ converge to the intrinsic distance between $X$ and $Y$ on $\Sigma$. This implies that the intrinsic diameters $D_n$ of the $\partial C_n$ converge to $D$.  Indeed, choosing $X,Y\in\Sigma$ at distance $D$ and points $X_n\to X$, $Y_n\to Y$ on $\partial C_n$ shows that $\liminf D_n\geq D$. Conversely, choose $X_n,Y_n\in\partial C_n$ at distance $D_n$; along a subsequence they converge to points of $\Sigma$, whose distance is at most $D$, so $\limsup D_n\leq D$. Surface area is continuous in the Hausdorff metric~\cite[Theorem~4.2.1 and (4.9)]{Schneider}, so the areas of the $\partial C_n$ converge to $A$. Passing to the limit gives the claim.
\end{proof}

\section{Relation with Shioya's estimate}\label{sec:shioya}

Shioya's proof~\cite{Shioya} starts from the same bound \eqref{eq:layer} at every point $p$ and integrates it over $p\in M$ before replacing $R_p$ by a constant. On the left this gives $A^2$. After integrating out the auxiliary variables, the curvature term on the right is
\[
  \frac12\int_M\!\int_MK(q)\bigl(R_p-d(p,q)\bigr)^2\,dA(q)\,dA(p),
\]
the integrand of $P_D$ with $R_p$ in place of $D$. Because $R_p$ varies with $p$, this term does not separate under Fubini into an expression in the areas of balls; Shioya therefore bounds $\pi R_p^2$ by $\pi D^2$ and shrinks the region of integration in the curvature term using $R_p\geq R$, where $R=\min_pR_p$ is the radius of $M$. Bishop--Gromov comparison then gives $A/D^2\leq\pi(1-\rho^4/3)$ with $\rho=R/D$, and combining this with Bishop's inequality $A/D^2\leq\pi\rho^2$ gives his constant.

The present argument reverses the order.  It first averages the potential $P_D$, defined with the fixed scale $D$, which separates under Fubini and is controlled by \eqref{eq:BG}; only then does it correct from $D$ to the eccentricity $R_p$ at the single selected point.

\section{Surfaces of positive Euler characteristic}\label{sec:projective}

The method applies to every closed surface of positive Euler characteristic, that is, to the sphere and to the real projective plane.

\begin{proposition}\label{prop:rp2}
  Let $M$ be a smooth closed Riemannian surface with nonnegative Gaussian curvature, Euler characteristic $\chi>0$, area $A$ and diameter $D$. Then
  \[
    \frac{A}{D^2}\leq\pi\Bigl(1-\frac{\chi}{6}\Bigr).
  \]
\end{proposition}

\begin{proof}[Sketch of proof]
  For $\chi=2$ this follows from Theorem~\ref{thm:smooth}, so let $\chi=1$. Then $M$ is homeomorphic to $\mathbb RP^2$, it is not orientable, and Gauss--Bonnet gives $\int_MK\,dA=2\pi\chi$. The integrated Fiala--Hartman inequality holds on $M$ as well: its proof, by the first variation of length of distance circles and Gauss--Bonnet for balls, is local and does not use orientability, and it is used in exactly this setting in Shioya's proof of his bound for $\mathbb RP^2$~\cite[p.~282]{Shioya}. A proper open connected subset of $M$ is a noncompact surface, so its Euler characteristic is again at most~$1$, and Proposition~\ref{prop:FH} and Lemma~\ref{lem:layer} hold. The rest of the proof of Theorem~\ref{thm:smooth} is unchanged, with $4\pi$ replaced by $2\pi\chi$: Lemma~\ref{lem:average} gives $P_D(p)\geq\tfrac{\pi\chi}{3}D^2$, so $a=\sqrt{\pi\chi/3}\,D$ and $b=\sqrt{2\pi\chi}$, and the cut-off $a/b=D/\sqrt6$ is unchanged. For $0\leq\delta\leq D/\sqrt6$,
  \[
    F'(\delta)=-2\pi D+\pi\chi\sqrt{2/3}\,D+2\pi(1-\chi)\delta<0,
  \]
  and $F(\delta)=\pi(D-\delta)^2$ is again decreasing for $D/\sqrt6\leq\delta\leq D$. Hence $A/D^2\leq\pi-\tfrac12\cdot\tfrac{\pi\chi}{3}=\pi(1-\chi/6)$.
\end{proof}

Proposition~\ref{prop:rp2} should be compared with Shioya's bound $(\sqrt{9+6\chi}-3)\pi/\chi$~\cite[Corollary~2]{Shioya}. For the real projective plane it gives
\[
  \frac{A}{D^2}\leq\frac{5\pi}{6}\approx2.618,
\]
against Shioya's $(\sqrt{15}-3)\pi\approx2.743$; the round projective plane has $A/D^2=8/\pi\approx2.546$.

\section*{Use of AI}

The main idea of this paper, the averaging argument, was found by the author without AI assistance, in the setting of plane triangulations of maximum degree~$6$. The translation of the argument to smooth surfaces, and the writing of the proof, were done with the help of an AI assistant (Claude, Anthropic). The author takes full responsibility for the content.

\end{document}